\documentclass[12pt]{amsart}
\usepackage{enumerate, amssymb}
\usepackage[colorlinks]{hyperref}
\usepackage{diagrams}
\diagramstyle[scriptlabels,height=8mm,width=8mm]

\def\bdi{\begin{diagram}}
\def\edi{\end{diagram}}

\theoremstyle{plain}

\newtheorem{thm}{Theorem}[section]
\newtheorem{cor}[thm]{Corollary}
\newtheorem{lem}{Lemma}
\newtheorem{prop}[thm]{Proposition}

\theoremstyle{definition}
\newtheorem{defi}[thm]{Definition}
\newtheorem{defis}[thm]{Definitions}
\newtheorem{conj}[thm]{Conjecture}
\newtheorem{conv}[thm]{Convention}
\newtheorem{nota}[thm]{Notation}
\newtheorem{rem}[thm]{Remark}
\newtheorem{rems}[thm]{Remarks}
\newtheorem{exa}[thm]{Example}
\newtheorem{exas}[thm]{Examples}
\newtheorem{sit}[thm]{}

\newcommand{\brem}{\begin{rem}}
\newcommand{\brems}{\begin{rems}}
\newcommand{\erem}{\end{rem}}
\newcommand{\erems}{\end{rems}}
\newcommand{\bexa}{\begin{exa}}
\newcommand{\bexas}{\begin{exas}}
\newcommand{\eexa}{\end{exa}}
\newcommand{\eexas}{\end{exas}}
\newcommand{\bdefi}{\begin{defi}}
\newcommand{\edefi}{\end{defi}}
\newcommand{\bdefis}{\begin{defis}}
\newcommand{\edefis}{\end{defis}}
\newcommand{\bcor}{\begin{cor}}
\newcommand{\ecor}{\end{cor}}
\newcommand{\blem}{\begin{lem}}
\newcommand{\elem}{\end{lem}}
\newcommand{\bconv}{\begin{conv}}
\newcommand{\econv}{\end{conv}}
\newcommand{\bconj}{\begin{conj}}
\newcommand{\econj}{\end{conj}}
\newcommand{\bprop}{\begin{prop}}
\newcommand{\eprop}{\end{prop}}
\newcommand{\bthm}{\begin{thm}}
\newcommand{\ethm}{\end{thm}}
\newcommand{\bnota}{\begin{nota}}
\newcommand{\enota}{\end{nota}}
\newcommand{\bsit}{\begin{sit}}
\newcommand{\esit}{\end{sit}}
\newcommand{\be}{\begin{eqnarray}}
\newcommand{\ee}{\end{eqnarray}}
\newcommand{\bproof}{\begin{proof}}
\newcommand{\eproof}{\end{proof}}
\def\ba{\begin{array}}
\def\ea{\end{array}}

\def\bnum{\begin{enumerate}}
\def\enum{\end{enumerate}}

\newcommand{\C}{{\mathbb C}}

\newcommand{\Z}{{\mathbb Z}}

\newcommand{\G}{{\Gamma}}

\newcommand{\nlin}{\unitlength1mm\begin{picture}(0,9.25)
                       \put(0,0.75){\line(0,1){8.5}}
                      \end{picture}}
\newcommand{\slin}{\unitlength1mm\begin{picture}(0,0)
                       \put(0,-0.75){\line(0,-1){8.5}}
                      \end{picture}}

\newcommand{\vlin}[1]{\hspace{0.75mm}\unitlength1mm\begin{picture}(#1,0)
                       \put(0,0){\line(1,0){#1}}
                      \end{picture}\hspace{0.75mm}\rule[-3mm]{0mm}{4mm}}

\def\llin{\vlin{11.5}}
\newcommand{\lin}{\vlin{8.5}}

\newcommand{\co}[1]{\unitlength1mm\begin{picture}(0,8)
    \put(0,0){\circle{1.5}}
    \put(0,3){\makebox(0,5)[b]{$#1$}}
                      \end{picture}}

\newcommand{\crl}[2]{\unitlength1mm\begin{picture}(0,8)
    \put(0,0){\circle{1.5}}
    \put(-5,0){\makebox(0,5)[b]{$#1$}}
   \put(5,0){\makebox(0,5)[b]{$#2$}}
      \end{picture}
      \rule[-7mm]{0mm}{7mm}}

\newcommand{\cshiftup}[2]{\unitlength1mm\begin{picture}(0,9.25)
                       \put(0,10){\crl{#1}{#2}}
                      \end{picture}}

\newcommand{\cshiftdown}[2]{\unitlength1mm\begin{picture}(0,9.25)
                       \put(0,-10){\crl{#1}{#2}}
                      \end{picture}}

\title[Corrigendum]{Corrigendum to our paper: \\Birational transformations of weighted
graphs}

\author{Hubert Flenner}
\address{Fakult\"at f\"ur Mathematik,
Ruhr Universit\"at Bochum, Geb.\ NA 2/72, Universit\"ats\-str.\
150, 44780 Bochum, Germany}
\email{Hubert.Flenner@rub.de}

\author{Shulim Kaliman}
\address{Department of Mathematics,
University of Miami, Coral Gables, FL  33124, U.S.A.}
\email{kaliman@math.miami.edu}

\author{Mikhail Zaidenberg}
\address{Universit\'e
Grenoble I, Institut Fourier, UMR 5582 CNRS-UJF, BP 74, 38402
St.\ Martin d'H\`eres c\'edex, France}
\email{zaidenbe@ujf-grenoble.fr}

\begin{document}

\begin{abstract}
We give a corrected version of Corollary 3.33 in \cite{FKZ1}.
\end{abstract}

\maketitle


In our paper \cite{FKZ1} on birational transformations of weighted
graph, we showed that every weighted graph $\G$ is birationally
equivalent  to a standard one.\footnote{See \S 2.1 and \S 2.4 in
{\em loc.cit.} for the definitions of birational transformations,
standard weighted graphs, and linear and circular segments. A
standard graph is supposed to be connected. } Concerning the
uniqueness of this standard form, we claimed in Corollary 3.33
that {\em a non-circular standard weighted graph is unique in its
birational equivalence class up to reversion of its linear
segments.}

However in this form the corollary is not true as was
pointed out to us by Karol Palka. The
correct version is as follows.
\medskip

\noindent {\bf Corollary 3.33.} {\em Every non-circular
standard weighted graph is unique in its birational equivalence
class up to reversions of its linear segments {\bf\em and changing
the weights of its branching vertices}. Similarly, a circular
standard graph is unique in its birational equivalence class up to
a cyclic permutation of its nonzero weights and reversion.}
\medskip

In the subsequent papers \cite{FKZ2}-\cite{FKZ4} we applied this
result only to linear or circular graphs. The latter graphs have
no branching vertices. Hence we used a correct version, and so
this does not affect any of our subsequent results.

To demonstrate that this correction is needed, following
a suggestion of Karol Palka let us consider the weighted tree

 
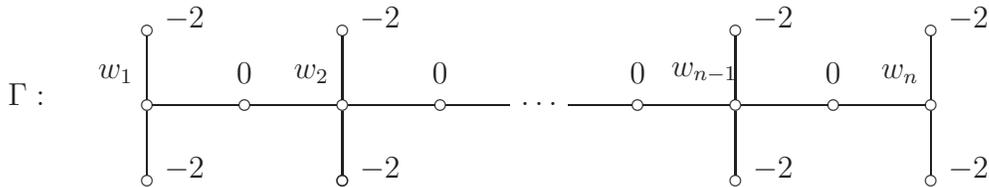
\vspace*{5mm}
$$
 \Gamma:\quad\qquad
 \cshiftup{}{-2}\nlin
 \slin\co{w_1\qquad}{}\cshiftdown{}{-2}
\llin\co{0}\llin\co{w_2\qquad}\slin \cshiftdown{}{}
\nlin\cshiftup{}{-2}\cshiftdown{}{-2}\llin\co{0}
\lin\ldots\lin\co{0}\llin\cshiftup{}{-2}\nlin
 \slin\co{w_{n-1}\quad\quad}{}\cshiftdown{}{-2}
\llin\co{0}\llin\co{w_n\qquad}\slin
\nlin\cshiftup{}{-2}\cshiftdown{}{-2} \quad
$$
\vspace*{10mm}

\noindent
Performing inner elementary transformations\footnote{See
\cite[\S 2.3]{FKZ1}.} at  zero vertices, one can change the
weights $w_1,\ldots, w_{n}$ of the branching vertices and so
replace them by the new weights $w'_1=\ldots w'_{n-1}=0$ and
$w_n'=\sum_{i=1}^n w_i$, without changing the birational class of
$\G$. Thus the correction above is indeed necessary.

It is possible to describe completely the possible weights
of branching points in the birational equivalence class of
a standard graph. For this we let $B$ denote the set
of branching points of $\G$.
Moreover let $\G_0$ be the subgraph of $\Gamma$, which
is the union of $B$ and all
linear segments of type $[[0_{2k+1}]]$ of $\G$ including
the edges that join these segments with $B$. For instance,
if $\G$ has no segments of type
$[[0_{2k+1}]]$ then $\G_0=B$.

With these notations the following more precise form
of Corollary 3.33 holds.
\medskip

\noindent
{\bf Corollary 3.33$'$.} {\em Let $\Gamma$ be a non-circular
standard weighted graph. Then $\G$  is unique in its birational
equivalence class up to reversion of its linear segments and the
change of weights of branching points described by the following
procedure:

\bnum\item[(1)] If a connected component $G$ of $\Gamma_0$
contains an end-vertex of $\Gamma$
\footnote{I.e., a vertex of $\Gamma$ of degree $1$.} then the
weights of points of $G\cap B$ can be chosen
arbitrarily;

\item[(2)] If $G$ does not contain any end-vertex of $\Gamma$ then the
weights of the points of $G\cap B$ can be chosen
arbitrarily modulo preservation of the sum of weights at these
branching points. \enum

\noindent In particular, we can change the weights of $\G$ in such
a way that for every connected component $G$ of $\G_0$ the new
weights of points of $G\cap B$  are all zero in case (1), and are
all zero with one exception in case (2), with the new weight in
this exceptional position equal to the sum of the weights of the
points of $G\cap B$.}
\medskip

For the proof we need a few preparations. Let us recall first
the following results.

\bnum\item[(i)] Any birational transformation of one standard
graph into another one is composed of a sequence of elementary
transformations (\cite[Theorem 3.1]{FKZ1}). In particular, it
identifies the branching vertices of both graphs, and transforms
linear segments of the first graph into linear segments of the
second one.

\item[(ii)] A standard linear segment is unique in its birational
equivalence class up to reversion (\cite[Corollary 3.32]{FKZ1}).
The reversion of a standard linear segment $L$ of $\G$ leaves
unchanged the weighted graph $\G\ominus L$ (\cite[Lemma
2.12]{FKZ1}).


\enum

\setcounter{thm}{0}

Let $\sigma: \G\rDotsto \G'$ be a birational transformation of
weighted graphs and let $L$ be a segment of $\G$. We call
$\sigma$ an $L$-{\em transformation} if it is composed of blowups
and blowdowns of $L$ and their subsequent total transforms. Thus
$\sigma$ induces a bijection $\G\ominus L\to \G'\ominus L'$ that
respects the weights of all vertices except possibly the branching
points that are adjacent to $L$ and $L'$, respectively. We
need the following lemma.

\setcounter{thm}{0}

\blem
(a) Every birational transformation $\sigma: \G\rDotsto \G'$
of standard graphs can be written as a composition of $L$-transformations,
where $L$ runs through the segments of $\G$.

(b) Every birational transformation $\sigma: \G \rDotsto \G'$
of standard graphs admits a domination
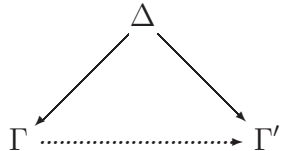
\bdi[small]
&&\Delta\\
&\ldTo & &\rdTo \\
\G&& \rDotsto &&\G' \edi such that $\Delta$, $\G$ and $\G'$ have
the same branching points. In particular, if $L$ is a segment of
$\G$ and $L'$ is the corresponding segment of $\G'$ then
the induced birational transformation $\sigma|L: L \rDotsto L'$
can be dominated by a linear graph. \elem

\bproof By (i) $\sigma$ can be written as a composition of
elementary transformations $\sigma_1\cdot\ldots\cdot \sigma_n$.
Every such elementary transformation takes place at some segment
$L$ of $\G$. Obviously elementary transformations taking place at
different segments commute. Thus if $\sigma_L$ is the product of
all $\sigma_i$ taking place at $L$ then $\sigma$ is the product of
all $\sigma_L$. This proves (a).

The first part of (b) follows from Lemma 3.8 in
\cite{FKZ1} while the second is immediate from the first one.
\eproof

\bproof[Proof of Corollary 3.33'.] Let $\sigma: \G \rDotsto \G'$
be a birational transformation of standard graphs.  By Lemma 1(b)
they can be dominated by a graph $\Delta$ such that $\Delta$, $\G$
and $\G'$ have the same branching points. Moreover by Lemma 1(a)
$\sigma$ can be decomposed as a product $\sigma =\prod_L\sigma_L$,
where $\sigma_L$ is an  $L$-transformation.

Let us first describe the effect of $\sigma_L$ on
the weights of the branching points of $\G$.

(a) If a branching point is not adjacent to $L$
then its weight remains unchanged under $\sigma_L$.

(b) Assume that $L= [[0_{2k}, w_1,\ldots, w_n]]$
with $n\ge 0$ and $w_i\le -2$ for all $i$.
According to Proposition 3.4 in \cite{FKZ1}
the induced birational transformation $\sigma|L=\sigma_L|L:L\rDotsto  L'$
is either the reversion or the identity. As observed in (ii),
in both cases $\sigma_L$ does not affect the weights of branching points.

(c) If $L=[[0_{2k+1}]]$ then according to Proposition 3.7
in \cite{FKZ1} $\sigma_L|L=\tau^s$ for some $s\in\Z$,
where $\tau$ denotes the left move (see Definition 3.5 in {\em loc.cit.}).
Comparing with Lemma 2.12(c) in {\em loc.cit.}, $\sigma_L$
will increase the weight of one of the branching points adjacent
to $L$ by $s$, while it decreases the weight of the other one
(if existent) by $s$.

In any case, for every connected component $G$ of $\G_0$
the sum of the weights of points in $G\cap B$ is invariant
under $\sigma_L$ and hence also under $\sigma$.

Finally, applying (c) repeatedly we can move the weight of any
vertex $v$ in $G\cap B$ to any other vertex $v'$ in $G\cap B$
along a path connecting them in $G$. To annihilate all the weights
but one, we proceed as follows. We choose a rooted subtree of $G$
containing all the vertices of $G\cap B$ (and the root is one of
them), and move recursively the weights to the root starting from
the vertices in $G\cap B$ with maximal distance from the root. In
the case that $G$ contains an end vertex of $\G$, choosing for the
root the vertex in $G\cap B$ closest to the end vertex, again by
(c) we finish the recursive procedure by annihilating the weight
of the root.\footnote{This procedure is transparent in the example
above. } This proves Corollary 3.33$'$.\eproof

\noindent {\bf Acknowledgments.} It is pleasure for us to thank Karol
Palka, who attracted our attention to the flaw in \cite{FKZ1}. We
also grateful to Daniel Daigle, who indicated in a letter to the
authors that in the new formulation, Corollary 3.33$'$ (for the
particular case where $\G$ is a forest) is essentially equivalent
to the classification result obtained in \cite{Dai2}. The case of
a linear graph is contained also in \cite{Dai1}. The relation to
the earlier classification due to Danilov and Gizatullin
\cite{DaGi} is explained in \cite{FKZ1}.

\end{document}